\documentclass[11pt, reqno]{amsart}
\usepackage{amsmath,amsfonts,amssymb,amsthm,color,cancel}
\usepackage{epsfig}
\usepackage[normalem]{ulem}
\usepackage{ulem}

\usepackage{bbm}
\usepackage{enumerate}
\usepackage{amstext}
\usepackage{mathrsfs}
\usepackage{xspace}
\usepackage{amsfonts}
\usepackage{amsmath}
\usepackage{amssymb}
\usepackage{amstext}
\usepackage{amsthm}
\usepackage{xspace}
\usepackage[active]{srcltx}
\usepackage{cite}

\newtheorem{theorem}{Theorem}[section]
\newtheorem{lemma}[theorem]{Lemma}

\newcommand{\N}{\mathbb{N}}
\newcommand{\Z}{\mathbb{Z}}
\newcommand{\R}{\mathbb{R}}
\newcommand{\Q}{\mathbb{Q}}
\newcommand{\T}{\mathbb{T}}

\newcommand{\Proof}{\textbf{Proof}\hspace{0.3cm}}
\newcommand{\End}{\ensuremath{\hfill{\Box}}\\}

\numberwithin{equation}{section}

\begin{document}




\title{QUANTITATIVE DESTRUCTION OF INVARIANT CIRCLES}

\author{Lin Wang}

\address{School of Mathematics and Statistics, Beijing Institute of Technology, Beijing 100081, China}

\email{lwang@bit.edu.cn}


\begin{abstract}
\noindent For
area-preserving twist maps on the annulus, we consider the problem on  quantitative destruction of invariant circles with a given frequency $\omega$ of an integrable system by a trigonometric
polynomial of degree $N$ perturbation $R_N$ with $\|R_N\|_{C^r}<\epsilon$. We obtain a relation among  $N$, $r$, $\epsilon$ and the arithmetic property of $\omega$, for which  the area-preserving map  admit no invariant circles with $\omega$.
\end{abstract}

\date{\today}
\maketitle



\section{Introduction and main result}
Area-preserving twist maps on the annulus first appeared in the work of Poincar\'{e}  on the three-body problem. They served as suitable prototypes for the study of a complicated Hamiltonian system.  As two dimensional discrete dynamical
models,  they
describe the behavior of area preserving surface diffeomorphisms in
the neighborhood of a generic elliptic periodic point. The study of such maps was initiated by Birkhoff in the 1920s. Since then, this class of maps have offered many
opportunities for the rigorous analysis of aspects of Hamiltonian systems.  As a highlight, Moser  proved the first
differentiable version of the KAM theorem in the context of twist maps.

The smoothness of the perturbation in the KAM theorem can be reduced. By the efforts of Moser,  Takens, R\"{u}ssman  and Herman   {\it et al.}, it was proved that  certain invariant circle with constant type frequency can
be persisted under arbitrarily small perturbations in the $C^3$
topology, where the invariant circle (also referred to as an essential curve) is an invariant curve that is not homotopic
to a point.

The works on non-existence of the invariant circles belong to converse KAM theory.
It was shown by
Herman  \cite{H2} that the invariant circle with a given frequency can be destroyed by $C^{3-\delta}$ arbitrarily small
$C^\infty$ perturbations.
Following the ideas and techniques
developed by Mather in the 1980s, a variational proof of Herman's result was provided in
\cite{W}. As a complement, it was considered in \cite{Wg}  for Gevrey-$\alpha$ ($\alpha>1$) systems to destroy the invariant circles with given frequencies. For Hamiltonian systems with multi-degrees of freedom, the corresponding results were obtained by \cite{CW} and \cite{P}. Moreover, it was obtained that all of the Lagrangian tori of  an integrable positive definite Hamiltonian system with $d$ ($d\geq 2$) degrees of freedom can be destroyed by an arbitrarily small  $C^\omega$ perturbation in the $C^{d-\delta}$ topology \cite{W3}.

For certain special frequencies, it was
obtained by Mather (resp. Forni) in \cite{M4} (resp. \cite{F}) that
the invariant circles with those frequencies can be destroyed by
small perturbations in finer topology respectively. More precisely,
Mather considered Liouvillian frequencies and the topology of the
perturbation induced by $C^\infty$ metric. Forni was concerned more about
more special frequencies which can be approximated  by rational
ones exponentially and the topology of the perturbation induced by
the supremum norm of real-analytic function.  Recently, Chen and Cheng \cite{CC} gave an open and dense property about the destruction of invariant
circles by using regular dependence of the Peierls barriers on
perturbations. Roughly speaking, there
is a balance between the arithmetic property of the frequency,
the regularity of the perturbation and its topology.

From physical point of view, it is more natural to consider real analytic perturbations, e.g. trigonometric
polynomials instead of $C^\infty$ ones.  Consider a completely integrable system with
the generating function
\[h_0(x,x')=\frac{1}{2}(x-x')^2, \quad x,x'\in \R.\]
Based on Poincar\'{e}' pioneering work, it is well known that if $\omega\in \Q$, then the invariant circles with frequency $\omega$ could be easily destroyed by an analytic
perturbation arbitrarily close to $0$ in the  topology induced by
the supremum norm of real-analytic function. Therefore it suffices to
consider the irrational $\omega$. Herman proved in \cite{H2} that the invariant circle with a given irrational frequency is unique.

 An
irrational number $\omega\in\R$ is called $\mu$-well approximable if
there exist infinitely many
integers $q_n\in\N$ such that
\begin{equation}\label{mu app}
|q_n\omega-p_n|<{q_n}^{-1-\mu},
\end{equation}
for some integer $p_n$.
It follows from Dirichlet approximation that any irrational number
is 0-well approximable.  $\omega$ is called a Liouvillian number if it
is $\mu$-well approximable for all $\mu>0$. Otherwise, it is called a Diophantine number.
 Moreover, Jarn\'{\i}k's theorem  shows that the set of $\mu$-well approximable numbers has Hausdorff dimension $\frac{2}{2+\mu}$.

Compared to the aforementioned results above, one can  ask the
following question from a quantitative of view: \begin{quote}
\textbf{Question:}
Given an irrational frequency $\omega$ and $0<\epsilon\ll 1$. Let $R_N(x)$ be a trigonometric
polynomial of degree $N$ which satisfies $\|R_N\|_{C^r}<\epsilon$. If the area-preserving map generated by $h_0(x,x')+R_N(x')$ admit no invariant circles with $\omega$, then what are the relation among $\epsilon$, $N$, $r$ and the arithmetic property of $\omega$?
\end{quote}

Based on the KAM result (\cite{H33}), we know that $r<4$ for a given badly approximable frequency $\omega$ if the perturbation is taken among $C^\infty$ functions. Nevertheless, it is not natural to expect the same upper bound of $r$ for a trigonometric
polynomial of degree $N$. Given $\mu$ and $\epsilon$, we are devoted to looking for a bigger $r$  and a smaller $N$. First of all, we obtain a relation between $r$ and $\mu$. Given any $u\in C^r(\mathbb{T})$. Let us recall the $C^r$-norm:
\[\|u\|_{C^r}:=\max_{|\alpha|\leq r}\max_{x\in \mathbb{T}}|D^{\alpha}u(x)|,\]
where $\alpha:=(\alpha_1,\ldots,\alpha_n)$, $|\alpha|:=\alpha_1+\cdots+\alpha_n$.
\begin{theorem}\label{liou}
Let $\omega$ be a $\mu$-well approximable frequency. Given $0<\epsilon\ll 1$, there exists   a trigonometric
polynomial $R_N(x)$
which satisfies $\|R_N\|_{C^r}<\epsilon$ with $r<3+\mu$, such that the area-preserving map generated by $h_0(x,x')+R_N(x')$ admit no invariant circles with frequency $\omega$.
\end{theorem}

Theorem \ref{liou} provides an analytic version of  Mather's result on non-existence of the invariant circle with a Liouvillian frequency in \cite{M4}. It implies that
both of small $C^\infty$ and trigonometric
polynomial perturbations play the same role in the $C^\infty$ topology for destroying invariant circles with given Liouvillian frequencies.

It is well known that the rigidity of the invariant circle with a $\mu$-well approximable  frequency  increases as $\mu$ decreases. In particular,
 $\omega$ is called badly approximable   if it
is exactly $0$-well approximable.  Based on Jarn\'{\i}k's theorem, the set of badly approximable numbers has Hausdorff dimension $1$ and Lebesgue measure $0$.

In order to find a smaller $N$, we provide a relation between $N$ and $\epsilon$ for the the invariant circle with a badly approximable  frequency.
\begin{theorem}\label{bada}
Let  $\omega$ be a  badly approximable frequency. Given $0<\epsilon\ll 1$ and $r\in [0,3)$, there exists $R_N(x)$ of degree \[N\leq C\epsilon^{-\frac{3}{2(3-r)}}\] with $\|R_N\|_{C^r}<\epsilon$ such that $h_0(x,x')+R_N(x')$ admit no invariant circles with frequency $\omega$.
\end{theorem}

It is worth noting that we can not conclude the optimality of $r$ and $N$ in Theorem \ref{liou} and  Theorem \ref{bada}. Note that for badly approximable frequencies, we obtain $r<3$, which is smaller than $r=4$ in the KAM result in \cite{H33}. To fill the gap, some further developments of (converse)
KAM theory are needed.
\section{ Preliminaries}

\subsection{Minimal configuration}
Let $F$ be a diffeomorphism of $\R^2$ denoted by $F(x,y)=(X(x,y),Y(x,y))$. Let $F$ satisfy:
\begin{itemize}
\item {\it Periodicity:} $F\circ T=T\circ F$ for the translation $T(x,y)=(x+1,y)$;
\item {\it Twist condition:} the map $\psi:(x,y)\mapsto(x,X(x,y))$ is a diffeomorphism of $\R^2$;
\item {\it Exact symplectic:} there exists a real valued function $h$ on $\R^2$ with $h(x+1,y)=h(x,y)$ such that
    \[YdX-ydx=dh.\]
\end{itemize}
Then $F$ induces a map on the cylinder denoted by $f$: $\T\times\R\mapsto \T\times\R$ ($\T=\R/\Z$). $f$ is called an exact
area-preserving  twist map. The function  $h$: $\R^2\rightarrow\R^2$ is called
a generating function of $F$, namely $F$
is generated by the following equations
\begin{equation*}
\begin{cases}
y=-\partial_1 h(x,x'),\\
y'=\partial_2 h(x,x'),
\end{cases}
\end{equation*}
where $F(x,y)=(x',y')$.

The function $F$ gives rise to a dynamical
system whose orbits are given by the images of points of $\R^2$
under the successive iterates of $F$. The orbit of the point
$(x_0,y_0)$ is the bi-infinite sequence
\[\{...,(x_{-k},y_{-k}),...,(x_{-1},y_{-1}),(x_0,y_0),(x_1,y_1),...,(x_k,y_k),...\},\]
where $(x_k,y_k)=F(x_{k-1},y_{k-1})$. The sequence
\[(...,x_{-k},...,x_{-1},x_0,x_1,...,x_k,...)\] denoted by $(x_i)_{i\in\Z}$ is called a
stationary configuration if it satisfies the identity
\[\partial_1 h(x_i,x_{i+1})+\partial_2 h(x_{i-1},x_i)=0,\ \text{for\ every\ }i\in\Z.\]
Given a sequence of points $(z_i,...,z_j)$, we can associate its
action
\[h(z_i,...,z_j)=\sum_{i\leq s<j}h(z_s,z_{s+1}).\] A configuration $(x_i)_{i\in\Z}$
is called minimal if for any $i<j\in \Z$, the segment
$(x_i,...,x_j)$ minimizes $h(z_i,...,z_j)$ among all segments
$(z_i,...,z_j)$ of the configuration  satisfying $z_i=x_i$ and
$z_j=x_j$. It is easy to see that every minimal configuration is a
stationary configuration. There is a visual way to describe  configurations. A configuration $(x_i)_{i\in\Z}$ is a function from $\Z$ to $\R$. One can interpolate this function linearly and obtain a piecewise affine function $\R\rightarrow\R$ denoted by $t\mapsto x_t$. The graph of this function is sometimes called the Aubry diagram of the configuration.
By \cite{B}, minimal configurations satisfy a
group of remarkable properties as follows:
\begin{itemize}
\item Two distinct minimal configurations seen as the Aubry diagrams cross at most once, which
is so called Aubry's crossing lemma.
\item For every minimal configuration $\bold{x}=(x_i)_{i\in\Z}$, the limit
\[\rho(\bold{x})=\lim_{n\rightarrow\infty}\frac{x_{i+n}-x_i}{n}\]
exists and doesn't depend on $i\in\Z$. $\rho(\bold{x})$ is called
the frequency of $\bold{x}$.
\item For every $\omega\in \R$, there exists a minimal configuration
with frequency $\omega$. Following the notations of \cite{B}, the
set of all minimal configurations with frequency $\omega$ is
denoted by $M_\omega^h$, which can be endowed with the topology
induced from the product topology on $\R^\Z$. If
$\bold{x}=(x_i)_{i\in\Z}$ is a minimal configuration, considering
the projection $pr:\ M_\omega^h\rightarrow\R$ defined by
$pr(\bold{x})=x_0$, we set $\mathcal {A}_\omega^h=pr(M_\omega^h)$.
\item If $\omega\in\Q$, say $\omega=p/q$ (in lowest terms), then it is convenient to define the rotation symbol to detect the structure of
$M_{p/q}^h$. If $\bold{x}$ is a minimal configuration with frequency $p/q $, then the rotation symbol $\sigma(\bold{x})$ of
$\bold{x}$ is defined as follows
\begin{equation*}
\sigma(\bold{x})=\left\{\begin{array}{ll}
\hspace{-0.4em}p/q+,&\text{if}\ x_{i+q}>x_i+p\ \text{for\ all\ }i,\\
\hspace{-0.4em}p/q,&\text{if}\ x_{i+q}=x_i+p\ \text{for\ all\ }i,\\
\hspace{-0.4em}p/q-,&\text{if}\ x_{i+q}<x_i+p\ \text{for\ all\ }i.\\
\end{array}\right.
\end{equation*}
 Moreover, we set
\begin{align*}
&M_{{p/q}^+}^h=\{\bold{x} \text{\  is a minimal configuration with
rotation symbol}\  p/q \text{\ or\ } p/q+\},\\
&M_{{p/q}^-}^h=\{\bold{x} \text{\ is a minimal configuration with
rotation symbol}\  p/q \text{\ or\ } p/q-\},
\end{align*}
then both $M_{{p/q}^+}^h$ and $M_{{p/q}^+}^h$ are totally ordered.
Namely, every two configurations in each of them (seen as Aubry diagrams) do not cross. We
denote $pr(M_{{p/q}^+}^h)$ and $pr(M_{{p/q}^-}^h)$ by $\mathcal
{A}_{{p/q}^+}^h$ and $\mathcal {A}_{{p/q}^-}^h$ respectively.
\item If $\omega\in\R\backslash\Q$ and $\bold{x}$ is a minimal
configuration with frequency $\omega$, then
$\sigma(\bold{x})=\omega$ and $M_\omega^h$ is totally ordered.
\item $\mathcal {A}_\omega^h$ is a closed subset of $\R$ for every rotation symbol
$\omega$.
\end{itemize}
\subsection{Peierls's barrier}
In \cite{M3}, Mather introduced the notion of Peierls's barrier and gave
a criterion of existence of invariant circle. Namely, the exact
area-preserving  twist map generated by $h$ admits an
invariant circle with frequency $\omega$ if and only if the
Peierls's barrier $P_\omega^h(\xi)$ vanishes identically for all
$\xi\in\R$. The Peierls's barrier is defined as follows:
\begin{itemize}
\item If $\xi\in \mathcal {A}_\omega^h$, we set $P_\omega^h(\xi)$=0.
\item If $\xi \not\in \mathcal {A}_\omega^h$, since $\mathcal {A}_\omega^h$ is a closed set in $\R$, then $\xi$ belongs to some
complementary interval $(\xi^-,\xi^+)$ of $\mathcal {A}_\omega^h$ in
$\R$. By the definition of $\mathcal {A}_\omega^h$, there exist
minimal configurations with rotation symbol $\omega$,
$\bold{x^-}=(x_i^-)_{i\in\Z}$ and $\bold{x^+}=(x_i^+)_{i\in\Z}$
satisfying $x_0^-=\xi^-$ and $x_0^+=\xi^+$. For every configuration
$\bold{x}=(x_i)_{i\in\Z}$ satisfying $x_i^-\leq x_i\leq x_i^+$, we
set
\[G_\omega(\bold{x})=\sum_I(h(x_i,x_{i+1})-h(x_i^-,x_{i+1}^-)),\]
where $I=\Z$, if $\omega$ is not a rational number, and $I=\{0,...,
q-1 \}$, if $\omega=p/q$. $P_\omega^h(\xi)$ is defined as the
minimum of $G_\omega(\bold{x})$ over the configurations $\bold{x}\in
\Pi=\prod_{i\in I}[x_i^-,x_i^+]$ satisfying $x_0=\xi$. Namely
\[P_\omega^h(\xi)=\min_{\bold{x}}\{G_\omega(\bold{x})|\bold{x}\in \Pi\ \text{and}\ \ x_0=\xi\}.\]
\end{itemize}
By \cite{M3}, $P_\omega^h(\xi)$ is a non-negative periodic function of
the variable $\xi\in\R$ with the modulus of continuity with respect
to $\omega$ and its modulus of continuity with respect to $\omega$ can be bounded from above. Due to the periodicity of $P_\omega^h(\xi)$ with
respect to $\xi$, we only need to consider it in the interval $[0,1]$.

For simplicity, we don't distinguish the constant $C$ in the following
different estimate formulas unless it is necessary. In the following sections, we will
prove Theorem \ref{liou} and Theorem \ref{bada}, which will be achieved by the detailed analysis on the regularity
of Peierls's barrier and an approximation from trigonometric
polynomials to $C^\infty$ functions.

\section{ Construction of the generating functions}
 We construct the perturbation of $h_0(x,x')$ as follows. The first one is
\begin{equation}\label{31}
u_n(x)=\frac{1}{n^a}(1-\cos(2\pi x) ),\quad x\in \R,\end{equation}
where $n\in \N$ and $a$ is a positive constant independent of $n$.

Let $\bar{h}_n(x,x')=h_0(x,x')+u_n(x')$, we have
\begin{lemma}\label{lowstep} Let $(x_i)_{i\in \Z}$ be a minimal
configuration of $\bar{h}_n$ with rotation symbol $0^+$, then
\[x_{k+1}-x_k\geq C(n^{-\frac{a}{2}}),\quad \text{for}\quad x_k\in \left[\frac{1}{4},\frac{3}{4}\right].\]\end{lemma}

The proof of Lemma \ref{lowstep} is similar to \cite[Lemma 4.1]{W}. For the sake of completeness, we will give it in Appendix A.

We construct  the second part of the perturbation in the following.
Let $p_N(x)$ be a trigonometric polynomial of degree $N$. By Hadamard's three-circle theorem (see \cite[Page 286-287]{F} for more details), one has
that for any $r>0$,
\begin{equation}
||p_N(x)||_r\leq e^{rN}||p_N(x)||,
\end{equation}
where $||p_N(x)||_r$ denotes the maximum of $|p_N(z)|$ in the strip
$S_r=\{z\in\mathbb{C}|\,|\text{Im}z|\leq r\}$ of width $2r$ in the
complex plane and $||p_N(x)||$ denotes the maximum of $|p_N(x)|$ on
the real line. Without loss of generality, we take $r=1$, namely
\begin{equation}
||p_N(x)||_1\leq e^{N}\max|p_N(x)|.
\end{equation}
Then, by the Cauchy estimates, for any fixed $s\geq 0$, we have
\begin{equation}\label{cachy}
||p_N(x)||_{C^s}\leq C_se^N\max|p_N(x)|,
\end{equation}
where $C_s$ is a constant depending on $s$ only.

Based on Lemma \ref{lowstep}, we need to construct a real analytic
function with a ``bump" in correspondence with the interval
$\Lambda_n$ satisfying
\begin{equation}\label{length}
\mathcal {L}(\Lambda_n)\sim n^{-\frac{a}{2}}\quad\text{and}\quad
\Lambda_n\subset\left[\frac{1}{4},\frac{3}{4}\right],
\end{equation}
where $\mathcal {L}(\Lambda_n)$ denotes the Lebesgue measure of
$\Lambda_n$ and $f\sim g$ means that $\frac{1}{C}g<f<C g$ holds for
some constant $C>0$.

The ``bump" will be accomplished by using Jackson's approximation
theorem (see \cite[Theorem 13.6, p115]{Z}). Let $\phi(x)$ be a
$k$-times differentiable periodic function on $\R$, then for every
$N\in\N$, there exists a trigonometric polynomial $p_N(x)$ of degree
$N$ such that
\[\max|p_N(x)-\phi(x)|\leq A_kN^{-k}||\phi(x)||_{C^k},\]
where $A_k$ is a constant depending  on $k\in\N$ only.

We take a $C^\infty$ bump function $\phi$ supported on the interval
$\Lambda_n$, whose maximum is equal to $2$. By (\ref{length}), the
length of $\Lambda_n$ is bounded by $Cn^{-\frac{a}{2}}$. Thus, one can
choose $\phi(x)$ such that
\begin{equation}
||\phi(x)||_{C^k}\sim
\left(n^{\frac{a}{2}}\right)^k=n^{\frac{ak}{2}},
\end{equation}
where $k$ is determined by (\ref{kk}) below.
Then, chose $N$ large enough to achieve
\begin{equation}\label{N}
\sigma:=A_k N^{-k}||\phi(x)||_{C^k}\ll 1,
\end{equation}
where $\sigma$ is determined by (\ref{sig}) below. By Jackson's
approximation theorem, we can construct a trigonometric polynomial
$p_N(x)$ of degree $N$ such that:
\begin{equation}
\left\{\begin{array}{ll}\hspace{-0.4em}\max p_N(x)\geq 1,&\text{attained on}\ \Lambda_n,\\
\hspace{-0.4em}|p_N(x)|\leq\sigma,&\text{on}\ [0,1]\backslash\Lambda_n.\\
\end{array}\right.
\end{equation}
By (\ref{N}), we have
\begin{equation}\label{nn}
N\sim \sigma^{-\frac{1}{k}}n^{\frac{a}{2}}.
\end{equation}
Finally, we consider the normalized trigonometric polynomial
\begin{equation}
\tilde{p}_{2N}(x)=e^{-2N}\left(\frac{p_N(x)}{\max p_N(x)}\right)^2.
\end{equation}
From (\ref{cachy}), $\tilde{p}_N(x)$ satisfies:
\begin{equation}
\left\{\begin{array}{ll}\hspace{-0.4em}\tilde{p}_{2N}(x)\geq 0,\\
\hspace{-0.4em}||\tilde{p}_{2N}(x)||_{C^s}\leq C,\\
\hspace{-0.4em}\max \tilde{p}_{2N}(x)=e^{-2N},&\text{attained on}\ \Lambda_n,\\
\hspace{-0.4em}|\tilde{p}_{2N}(x)|\leq \sigma^2e^{-2N},&\text{on}\ [0,1]\backslash\Lambda_n.\\
\end{array}\right.
\end{equation}
Based on preparations above, we can construct the second part of the
perturbation as follow
\begin{equation}
v_n(x)=u_n(x)\tilde{p}_{2N}(x)=\frac{1}{n^a}(1-\cos 2\pi
x)\tilde{p}_{2N}(x).
\end{equation}
It is easy to see $v_n$ satisfies the following properties:
\begin{equation}\label{vn}
\left\{\begin{array}
{ll}\hspace{-0.4em}v_n(x)\geq 0,&\\
\hspace{-0.4em}||v_n(x)||_{C^s}\leq Cn^{-a},&\\
\hspace{-0.4em}\max v_n(x)\geq e^{-2N}n^{-a},&\text{attained on}\ \Lambda_n,\\
\hspace{-0.4em}|v_n(x)|\leq C\sigma^2e^{-2N}n^{-a},&\text{on}\
[0,1]\backslash\Lambda_n.
\end{array}\right.
\end{equation}

So far, we complete the construction of the generating function of
the nearly integrable system,
\begin{equation}\label{h}
h_n(x,x')=h_0(x,x')+u_n(x')+v_n(x'),
\end{equation}
where $n\in\N$.

\section{ Proof of Theorem \ref{liou}}
 First, we prove the
non-existence of invariant circles with a small enough frequency. More precisely, we have the following Lemma:
\begin{lemma}\label{MR} For $\omega\in\R\backslash\Q$ and $n$ large enough,
the exact area-preserving  twist map generated by $h_n$
admits no invariant circle with the frequency satisfying
\[|\omega|<n^{-a-\delta}, \] where $\delta$ is a small positive constant independent of $n$.
\end{lemma}
\Proof First of all, we estimate the lower bound of $P_{0^+}^{h_n}$.
 Let $(\xi_i)_{i\in \Z}$
be a minimal configuration of $h_n$ defined by (\ref{h}) with
rotation symbol $0^+$ satisfying $\xi_0=\eta$, where $\eta$
satisfies $v_n(\eta)=\max v_n(x)$ and let $(x_i)_{i\in \Z}$ be the
minimal configuration of
$\bar{h}_n(x_i,x_{i+1})=h_0(x_i,x_{i+1})+u_n(x_{i+1})$ with rotation
symbol $0^+$, then
\begin{align*}
\sum_{i\in
\Z}(h_n(&\xi_i,\xi_{i+1})-h_n(\xi_i^-,\xi_{i+1}^-))\\
&\geq v_n(\eta)+\sum_{i\in \Z}\bar{h}_n(\xi_i,\xi_{i+1})-\sum_{i\in
\Z}h_n(\xi_i^-,\xi_{i+1}^-),\\
&\geq v_n(\eta)+\sum_{i\in \Z}\bar{h}_n(x_i,x_{i+1})-\sum_{i\in
\Z}h_n(x_i,x_{i+1}),\\
&=v_n(\eta)-\sum_{i\in \Z}v_n(x_{i+1}).\\
\end{align*}
By \cite[Page 208, (4.2)]{M4}, there holds
\[P_{0^+}^{h_n}(\eta)=\sum_{i\in
\Z}(h_n(\xi_i,\xi_{i+1})-h_n(\xi_i^-,\xi_{i+1}^-)).\]
Therefore, we have shown:
\begin{align*}
P_{0^+}^{h_n}(\eta)\geq
v_n(\eta)-\sum_{i\in \Z}v_n(x_{i+1}).
\end{align*}
By (\ref{vn}), we have
\begin{equation*}
v_n(\eta)\geq e^{-2N}n^{-a}.
\end{equation*}
It follows from (\ref{uwith}) that
\begin{equation*}
\sum_{i\in \Z}v_n(x_{i+1})\leq
\sigma^2e^{-2N}\sum_{i\in\Z}u_n(x_{i+1})\leq\sigma^2e^{-2N}\sum_{i\in\Z}\frac{1}{4}(x_{i+1}-x_{i-1})^2\leq\sigma^2e^{-2N}.
\end{equation*}
Hence,
\begin{equation*}
P_{0^+}^{h_n}(\eta)\geq e^{-2N}(n^{-a}-\sigma^2),
\end{equation*}
we choose then $\sigma$ (consequently $N$) in such a way that
\[\sigma^2=\frac{1}{4}n^{-a},\]
which means
\begin{equation}\label{sig}
\sigma=\frac{1}{2}n^{-\frac{a}{2}}.
\end{equation}
By (\ref{nn}), it follows that
\begin{equation}\label{kef}
N\sim n^{\frac{a}{2}+\frac{a}{2k}},
\end{equation}
from which we have
\begin{equation}\label{lowb}
P_{0^+}^{h_n}(\eta)\geq \frac{3}{4}
n^{-a}\exp\left(-Cn^{\frac{a}{2}+\frac{a}{2k}}\right).
\end{equation}

 Secondly, following  a similar argument as in \cite{W}, we have
\begin{equation}\label{app}
|P_{\omega}^{h_n}(\xi)-P_{0^+}^{h_n}(\xi)|\leq
C\exp\left(-2n^{\frac{a}{2}+\frac{\delta}{2}}\right).\end{equation}where
$\xi\in \Lambda_n$ and $\delta$ is a small positive constant
independent of $n$. Here, $\Lambda_n$ is  the same as
in (\ref{length}). For the sake of completeness, we will prove (\ref{app}) in Appendix B.

Based on the preparations above, it is easy to prove Lemma \ref{MR}.
We assume that there exists an invariant circle with frequency
$0<\omega<n^{-a-\delta}$ for $h_n$, then $P_\omega^{h_n}(\xi)\equiv
0$ for every $\xi\in \R$. By (\ref{app}), we have
\begin{equation}\label{pp}
|P_{0^+}^{h_n}(\xi)|\leq
C\exp\left(-2n^{\frac{a}{2}+\frac{\delta}{2}}\right),\quad\text{for}\quad
\xi\in \Lambda_n.
\end{equation}

On the other hand, $(\ref{lowb})$ implies that there exists a point
$\eta\in \Lambda_n$ such that
\[P_{0^+}^{h_n}(\eta)\geq \frac{3}{4}n^{-a}\exp\left(-Cn^{\frac{a}{2}+\frac{a}{2k}}\right).\]Hence, we have
\begin{equation}\label{cotr}
n^{-a}\exp\left(-Cn^{\frac{a}{2}+\frac{a}{2k}}\right)\leq
C\exp\left(-2n^{\frac{a}{2}+\frac{\delta}{2}}\right).
\end{equation}
To achieve the contradiction, it suffices to take
\begin{equation}\label{kk}
k>\frac{a}{\delta},
\end{equation}
which implies
\[\frac{a}{2k}<\frac{\delta}{2}.\]
Note that $C$ is independent of $n$. Hence, for $n$ large enough
\begin{equation*}
n^{-a}\exp\left(-Cn^{\frac{a}{2}+\frac{a}{2k}}\right)>
C\exp\left(-2n^{\frac{a}{2}+\frac{\delta}{2}}\right),
\end{equation*}
which contradicts (\ref{cotr}).
 Therefore,
there exists no invariant circle with frequency
$0<\omega<n^{-a-\delta}$.

For $-n^{-a-\delta}<\omega<0$, by comparing $P_{\omega}^{h_n}(\xi)$
with $P_{0^-}^{h_n}(\xi)$, the proof is similar. We omit the
details. This completes the proof of Lemma \ref{MR}. \End

The case with a given irrational frequency can be easily
reduced to the one with a small enough frequency. More
precisely,

\begin{lemma}\label{Herm} Let $h_P$ be a generating function as follows
\[h_P(x,x')=h_0(x,x')+P(x'),\] where $P$ is a periodic
function of periodic $1$. Let $Q(x)=q^{-2}P(qx),q\in \N$, then the
exact area-preserving  twist map generated by
$h_Q(x,x')=h_0(x,x')+Q(x')$ admits an invariant circle with frequency $\omega \in \R\backslash \Q$ if and only if the exact
area-preserving  twist map generated by $h_P$ admits an
invariant circle with frequency $q\omega-p, p\in \Z$.
\end{lemma}

We omit the proof and for more details, see \cite{H2}. For the sake of
simplicity of notations, we denote $Q_{q_n}$ by $Q_n$ and the same
for $u_{q_n}, v_{q_n}$ and $h_{q_n}$. Let
\[Q_n(x)={q_n}^{-2}(u_n(q_nx)+v_n(q_nx)),\] where $(q_n)_{n\in \N}$
is a sequence satisfying (\ref{mu app})
\begin{equation}\label{diri}
|q_n\omega-p_n|<\frac{1}{q^{1+\mu}_n}, \end{equation} where $p_n\in \Z$ and
$q_n\in \N$. Since $\omega\in\R\backslash\Q$, we have
$q_n\rightarrow\infty$ as $n\rightarrow\infty$. Let
$\tilde{h}_n(x,x')=h_0(x,x')+Q_n(x')$, we prove Theorem \ref{liou}
for $(\tilde{h}_n)_{n\in\N}$ as follow:

  From the constructions of $u_n$ and $v_n$,
it follows that
\begin{equation}\label{hnmin}
\begin{split}
 ||\tilde{h}_n&(x,x')-h_0(x,x')||_{C^r}\\
 &=||Q_n(x')||_{C^r},\\
&\leq{q_n}^{-2}(||u_n(q_nx')||_{C^r}+||v_n(q_nx')||_{C^r}),\\
&\leq{q_n}^{-2}({q_n}^{-a}(2\pi)^r{q_n}^r+C_1{q_n}^{-a}{q_n}^r),\\
&\leq C_2{q_n}^{r-a-2},
\end{split}
\end{equation}
where $C_1, C_2$ are positive constants  depending on $r$ only.

Hence, it is enough to make $r-a-2<0$.
Based on Lemma \ref{MR} and the Dirichlet approximation
(\ref{diri}), it suffices to take $a:=1+\mu-\delta$. Moreover, we choose
\[r<3+\mu,\quad \delta:=\frac{3+\mu-r}{2}.\]
Then
\[r-a-2=\frac{r-(3+\mu)}{2}<0.\]

This completes the
proof of Theorem \ref{liou}.

\section{ Proof of Theorem \ref{bada}}

 Given  a badly approximable frequency $\omega$. Let $p_n\in \Z$ and
$q_n\in \N$ satisfy (\ref{diri}). Let $a\in (0,1)$ and
\[R_{N'}(x):=\frac{1}{q_n^{2+a}}\left(1-\cos(2\pi q_n x)\right)\left(1+\tilde{p}_{2N}(q_n x)\right).\]
It is clear to see that $R_{N'}(x)$ is  a trigonometric
polynomial of degree $N'=(2N+1)q_n$. Next, we give an estimate on the upper bound of $N'$ in terms of $\epsilon$ and $r$.

 By (\ref{hnmin}), we have
\[\|R_{N'}(x)\|_{C^r}\leq C_1\left(\frac{1}{q_n}\right)^{2+a-r},\]
where $r\in [0,2+a)$. In order to achieve $\|R_{N'}(x)\|_{C^r}<\epsilon$, it suffices to require
\begin{equation}\label{q1}
q_n>C_2\epsilon^{-\frac{1}{2+a-r}}.
\end{equation}
In order to get a smaller $N'$, we assume
\[q_{n-1}\leq C_2\epsilon^{-\frac{1}{2+a-r}}.\]
It is clear that $\omega$ is  badly approximable  if and only if it is constant-type from continued fraction expansion point of view. Namely, there exists $K:=K(\omega)$ such that $|a_n|\leq K(\omega)$, where $a_n$ denotes the $n$th partial quotient of $\omega$.
By virtue of \cite[Lemma 5F]{S},  we have
\[q_n\leq C_3q_{n-1},\]
where $C_3$ is a positive constant independent of $n$. It follows that
\begin{equation}\label{q122}
q_n<C_4\epsilon^{-\frac{1}{2+a-r}}.
\end{equation}
By (\ref{kef}), to   ensure that the area-preserving map generated by $h_0(x,x')+R_{N'}(x')$ admit no invariant circles with frequency $\omega$, we only need
\[N\sim q_n^{\frac{a}{2}+\frac{a}{2k}},\]
where $k$ satisfies (\ref{kk}). Combining with (\ref{q1}), we have
\[N>C_5\epsilon^{-\frac{a(k+1)}{2k(2+a-r)}}.\]
Moreover, it gives rise to
\begin{equation}\label{q2}
N'=(2N+1)q_n>C_6\epsilon^{-\frac{a(k+1)}{2k(2+a-r)}-\frac{1}{2+a-r}}+C_2\epsilon^{-\frac{1}{2+a-r}}.
\end{equation}
We denote
\[\gamma:=\frac{a(k+1)}{2k(2+a-r)}+\frac{1}{2+a-r}=\frac{a(k+1)+2k}{2k(2+a-r)}.\]
For each $a\in (0,1)$ and $r\in [0,2+a)$, a direct calculation implies
\begin{equation}\label{gaaa}
\gamma=\frac{3}{2(3-r)}+\frac{(1-a)r}{2(2+a-r)(3-r)}+\frac{a}{2(2+a-r)}\frac{1}{k}.
\end{equation}

 Let $\delta:=1-a$. By (\ref{gaaa}), $\gamma$ can be reformulated as
 \begin{equation*}
\gamma=\frac{3}{2(3-r)}+\frac{\delta r}{2(3-r-\delta)(3-r)}+\frac{1-\delta}{2(3-r-\delta)}\frac{1}{k}.
\end{equation*}
  In terms of the Dirichlet approximation
(\ref{diri}), $\delta\in (0,1)$. We pursue the largest value of $r$ and the smallest value of $\gamma$. Thus, we focus on the case with $0<3-r\ll 1$. Moreover, one can take for $m\in \mathbb{N}_+$,
\[\delta=\frac{3-r}{2m},\]
It follows that
\begin{equation*}
\gamma=\frac{3}{2(3-r)}+\frac{r}{2(2m-1)(3-r)}+\frac{2m-3+r}{2(2m-1)(3-r)^2}\frac{1}{k}.
\end{equation*}
Given $r$ with $0<3-r\ll 1$ and  $0<\epsilon\ll 1$. One can take $k$ and $m$ large enough such that
\[\gamma<\frac{3}{2(3-r)}+\epsilon.\]
Note that $\epsilon^{\epsilon}\rightarrow 1$ as $\epsilon\rightarrow 0$. In view of (\ref{q122}) and (\ref{q2}), one can find a a trigonometric
polynomial of degree
\[N'\leq C\epsilon^{-\frac{3}{2(3-r)}},\]
such that the area-preserving map generated by $h_0(x,x')+R_{N'}(x')$ admit no invariant circles with frequency $\omega$, where $C$ is a positive constant independent of $\epsilon$.

This completes the
proof of Theorem \ref{bada}.

\appendix
\section{Proof of Lemma \ref{lowstep}}

 Without loss of generality, we assume $x_i\in [0,1]$ for all
$i\in\Z$. By Aubry's crossing lemma, we have
\[0<...<x_{i-1}<x_i<x_{i+1}<...<1.\]
Let $x_k\in \left[\frac{1}{4},\frac{3}{4}\right]$. We consider the configuration
$(\xi_i)_{i\in \Z}$ defined by
\begin{equation*}
\xi_i= \left\{\begin{array}{ll}\hspace{-0.4em}x_i,& i<k,\\
\hspace{-0.4em}x_{i+1},& i\geq k.\\
\end{array}\right.
\end{equation*}
Since $(x_i)_{i\in \Z}$ is minimal, we have
\[\sum_{i\in \Z}\bar{h}_n(\xi_i,\xi_{i+1})-\sum_{i\in \Z}\bar{h}_n(x_i,x_{i+1})\geq 0.\]
By the definitions of $\bar{h}_n$ and $(\xi_i)_{i\in\Z}$, we have
\begin{align*}
0&\leq\sum_{i\in \Z}\bar{h}_n(\xi_i,\xi_{i+1})-\sum_{i\in
\Z}\bar{h}_n(x_i,x_{i+1})\\
&=\bar{h}_n(x_{i-1},x_{i+1})-\bar{h}_n(x_{i-1},x_{i})-\bar{h}_n(x_{i},x_{i+1})\\
&=(x_{i+1}-x_i)(x_i-x_{i-1})-u_n(x_i).
\end{align*}
Moreover,\[u_n(x_i)\leq(x_{i+1}-x_i)(x_i-x_{i-1})\leq\frac{1}{4}(x_{i+1}-x_{i-1})^2.\]
Therefore,
\begin{equation}\label{uwith}
x_{i+1}-x_{i-1}\geq 2\sqrt{u_n(x_i)}. \end{equation} For $x_k\in
[\frac{1}{4},\frac{3}{4}]$, $u_n(x_k)\geq n^{-a}$, hence,
\begin{equation}\label{ls} x_{k+1}-x_{k-1}\geq 2n^{-\frac{a}{2}}.
\end{equation}

Since $(x_i)_{i\in\Z}$ is a stationary configuration, we have
\begin{align*}
x_{i+1}-x_i&=-\partial_1\bar{h}_n(x_i,x_{i+1}),\\
&=\partial_2\bar{h}_n(x_{i-1},x_i),\\
&=x_i-x_{i-1}+u_n'(x_i).
\end{align*}
Since $u_n'(x)=\frac{2\pi}{n^a}\sin(2\pi x)$, it follows from
$(\ref{ls})$ that
\[x_{k+1}-x_k\geq C(n^{-\frac{a}{2}}),\quad x_k\in \left[\frac{1}{4},\frac{3}{4}\right]
.\] The proof of Lemma \ref{lowstep} is completed.\End

\section{Proof of  (\ref{app})}
Following a similar argument as in \cite{F,M3,W}, we give some details of the proof of (\ref{app}). In \cite{W}, we obtain that
for certain $g_n$ and $\xi$ and every irrational rotation symbol $\omega$ satisfying $0<\omega<n^{-\frac{a}{2}-\delta}$,
\begin{equation}\label{B1}
|P_{\omega}^{g_n}(\xi)-P_{0^+}^{g_n}(\xi)|\leq
C\exp\left(-2n^{\frac{\delta}{2}}\right).\end{equation}
where $\delta$ is a small positive constant
independent of $n$. Here we aim to prove that
for every irrational rotation symbol $\omega$ satisfying $0<\omega<n^{-a-\delta}$,
\begin{equation}\label{B2}
|P_{\omega}^{h_n}(\xi)-P_{0^+}^{h_n}(\xi)|\leq
C\exp\left(-2n^{\frac{a}{2}+\frac{\delta}{2}}\right).\end{equation}
where
$\xi\in \Lambda_n$ and $\delta$ is a small positive constant
independent of $n$, and  $\Lambda_n$ is given by (\ref{length}).

To achieve (\ref{B1}), the key part is  to verify that each of the intervals $[0,\exp(-n^{\frac{\delta}{2}})]$ and $[1-\exp(-n^{\frac{\delta}{2}}),1]$ contains a large number of points of the minimal configuration of $g_n$ for $n$ large enough.

Correspondingly, in order to prove  (\ref{B2}), it suffices to show that each of the intervals $[0,\exp(-n^{\frac{a}{2}+\frac{\delta}{2}})]$ and $[1-\exp(-n^{\frac{a}{2}+\frac{\delta}{2}}),1]$ contains a large number of points of the minimal configuration of  $h_n$ for $n$ large enough.
Based on that, (\ref{app}) can be obtained by a standard argument (see \cite{W}). In the rest of the appendix, we are devoted to proving the following lemma.
\begin{lemma}\label{lbspaceA} Let $(x_i)_{i\in \Z}$ be a minimal
configuration of $h_n$  with rotation symbol
$0<\omega<n^{-a-\delta}$, then there exist $j^-,j^+\in\Z$
such that \begin{align*}&0<x_{j^--1}<x_{j^-}<x_{j^-+1}\leq
\exp(-n^{\frac{a}{2}+\frac{\delta}{2}}),\\
& 1-\exp(-n^{\frac{a}{2}+\frac{\delta}{2}})\leq
x_{j^+-1}<x_{j^+}<x_{j^++1}<1.\end{align*}\end{lemma}

To prove Lemma \ref{lbspaceA}, we need to do some preliminary work. First of
all, we count the number of the elements of a minimal configuration
$(x_i)_{i\in\Z}$ with arbitrary rotation symbol $\omega$ in a given
interval. With the method of \cite{F}, we can conclude the following
lemma.

\begin{lemma}\label{count 0} Let $(x_i)_{i\in\Z}$ be a minimal configuration
of $h_n$ with rotation symbol $\omega>0$,
$J_n=\left[\exp(-n^{\frac{a}{2}+\frac{\delta}{2}}),\frac{1}{2}\right]$
and $\Sigma_n=\{i\in \Z|\,x_i\in J_n\}$, then \[\sharp\Sigma_n\leq
Cn^{a+\frac{\delta}{2}},\] where $\sharp\Sigma_n$ denotes
the number of elements in $\Sigma_n$ and  $\delta$ is a small
positive constant independent of $n$. \end{lemma}
\Proof Let $x^-=\exp\left(-n^{\frac{a}{2}+\frac{\delta}{2}}\right),
x^+=\frac{1}{2}$ and
$\sigma=\left(\frac{x^+}{x^-}\right)^{\frac{1}{M}}$, hence,
\[\ln \sigma=\frac{\ln(x^+)-\ln(x^-)}{M}.\] We choose $M\in \N$
such that $1\leq\ln\sigma\leq 2$, then
$M\sim n^{\frac{a}{2}+\frac{\delta}{2}}$.

We consider the partition of the interval $J_n=[x^-,x^+]$ into the
subintervals $J_n^k=[\sigma^k x^-,\sigma^{k+1}x^-]$ where $0\leq
k<M$. Hence, $J_n=\cup_{k=0}^{M-1}J_n^k$. We set
$S_k=\{i\in\Sigma_n|(x_{i-1},x_{i+1})\subset J_n^k\}$ and
$m_k=\sharp S_k$.

By the similar deduction as the one in Lemma \ref{lowstep}, we have
\[x_{i+1}-x_{i-1}\geq 2\sqrt{u_n(x_i)+v_n(x_i)}\geq Cn^{-\frac{a}{2}}x_i,\quad\text{for}\quad x_i\in \left[0,\frac{1}{2}\right].\]
For simplicity of notation, we denote $Cn^{-\frac{a}{2}}$ by
$\alpha_n$.

If  there exists $k$ such that $i\in S_k$ for $(x_i)_{i\in \Z}$,
then $x_{i+1}-x_{i-1}\geq\alpha_n\sigma^kx^-$,
moreover,\[m_k\alpha_n\sigma^kx^-\leq 2\mathcal {L}(J_n^k)=
2(\sigma-1)\sigma^kx^-,\]where $\mathcal {L}(J_n^k)$ denotes the
length of the interval of $J_n^k$. Hence $m_k\leq 2(\sigma
-1)\alpha_n^{-1}$.

On the other hand, if $i\in \Sigma_n\backslash \cup_{k=0}^{M-1}S_k$
, then there exists $l$ satisfying $0\leq l<M$ such that
\[x_{i-1}<\sigma^lx^-<x_{i+1}.\]Hence,
\[\sharp\{i\in\Sigma_n|i\not\in S_k \ \text{for\ any}\ k\}\leq 2M.
\] Therefore,
\begin{align*}
\sharp(\Sigma_n)\leq 2M(\sigma-1)\alpha_n^{-1}+2M.
\end{align*}Since $1\leq\ln\sigma\leq 2$ and $M\sim n^{\frac{a}{2}+\frac{\delta}{2}}$, then we have
\[\sharp\Sigma_n\leq
Cn^{a+\frac{\delta}{2}}.\] The proof of Lemma \ref{count
0} is completed.\End

Let $(x_i)_{i\in\Z}$ be a minimal configuration of $h_n$ with rotation symbol $\omega>0$, An argument as similar
as the one in Lemma \ref{count 0} implies that
\[\sharp\left\{i\in\Z\bigg|\,x_i\in\left[\exp\left(-n^{\frac{a}{2}+\frac{\delta}{2}}\right),1-\exp\left(-n^{\frac{a}{2}+\frac{\delta}{2}}\right)\right]\right\}\leq Cn^{a+\frac{\delta}{2}}.\]

It is easy to count the number of the elements of a minimal
configuration with irrational rotation symbol. More precisely, we
have the following lemma.
 \begin{lemma}\label{count w}Let
$(x_i)_{i\in \Z}$ be a minimal configuration with frequency
$\omega\in \R\backslash\Q$. Then for every interval $I_k$ of length
$k$, $k\in \N$,
\[\frac{k}{\omega}-1\leq\sharp\{i\in \Z|x_i\in
I_k\}\leq\frac{k}{\omega}+1.\] \end{lemma}
\Proof For every minimal configuration $(x_i)_{i\in\Z}$ with frequency $\omega$, there exists an orientation-preserving circle
homeomorphism $\phi$ such that $\rho(\Phi)=\omega$, where
$\Phi:\R\rightarrow\R$ denotes a lift of $\phi$. Since $\omega\in
\R\backslash\Q$, thanks to \cite{H1}, $\phi$ has a unique invariant
probability measure $\mu$ on $\T$ such that $\int_x^{\Phi(x)}
d\mu=\omega$ for every $x\in\R$. We denote $\int_x^{\Phi(x)}
d\mu$ by $\mu([x,\Phi(x)])$. In particular,
\[\mu([x_i,x_{i+1}])=\omega,\quad \text{for\ every\ }i\in\Z.\]
From $\mu(I_k)=k$, it follow that
\begin{align*}
&\omega\cdot(\sharp\{i\in \Z|x_i\in I_k\}-1)\leq k,\\
&\omega\cdot(\sharp\{i\in \Z|x_i\in I_k\}+1)\geq k,
 \end{align*}
which completes the proof of Lemma \ref{count w}.\End

Based on Lemma \ref{count 0} and Lemma \ref{count w}, if
$0<\omega<n^{-a-\delta}$ and $\omega$ is irrational, then for $n$ large enough,
\begin{equation}\label{pnew}
\sharp\{i\in \Z|x_i\in I_1\}\geq \frac{1}{\omega}-1\geq
C_1n^{a+\delta}>C_2n^{a+\frac{\delta}{2}},
\end{equation}
where $I_1$ denotes the closed interval of length $1$.

Based on two counting lemmas above, it is easy to prove Lemma \ref{lbspaceA}.
 By contradiction, we assume that there exist at most two points
of $(x_i)_{i\in\Z}$ in $[0, \exp(-n^{\frac{a}{2}+\frac{\delta}{2}})]$, say $x_m$
and $x_{m+1}$. It follows that $x_{m-1}<0$ and
$x_{m+2}>\exp(-n^{\frac{a}{2}+\frac{\delta}{2}})$. Hence, among the intervals
$[x_{m-1}, x_m]$, $[x_m, x_{m+1}]$ and $[x_{m+1}, x_{m+2}]$, there
exists at least one such that its length is not less than
$\frac{1}{3}\exp(-n^{\frac{a}{2}+\frac{\delta}{2}})$. Without loss of
generality, say $[x_{m+1}, x_{m+2}]$.

Since $(x_i)_{i\in\Z}$ is a stationary configuration, we have
\[x_{m+2}-x_{m+1}=x_{m+1}-x_m+u'_n(x_{m+1}),\] where $u'_n(x_{m+1})=\frac{2\pi}{n^a}\sin(2\pi
x_{m+1})$. From $x_{m+1}\in [-\exp(-n^{\frac{a}{2}+\frac{\delta}{2}}),
\exp(-n^{\frac{a}{2}+\frac{\delta}{2}})]$, it follows that
\[|u'_n(x_{m+1})|\leq
Cn^{-a}\exp(-n^{\frac{a}{2}+\frac{\delta}{2}}),\]which implies there exists $K$
independent of $n$ such that $[-\exp(-n^{\frac{a}{2}+\frac{\delta}{2}}),
\exp(-n^{\frac{a}{2}+\frac{\delta}{2}})]$ contains at most $K$ points of
$(x_i)_{i\in \Z}$.

On the other hand, by (\ref{pnew}), we have that for $n$ large
enough, the number of points of $(x_i)_{i\in\Z}$ in
$[-\exp(-n^{\frac{a}{2}+\frac{\delta}{2}}), \exp(-n^{\frac{a}{2}+\frac{\delta}{2}})]$ is
also large enough, which is a contradiction. Therefore, there exists
$j^-\in\Z$ such that
\begin{align*} 0\leq
x_{j^--1}<x_{j^-}<x_{j^-+1}<\exp(-n^{\frac{a}{2}+\frac{\delta}{2}}).\end{align*}

Similarly, there exists $j^+\in\Z$ such that
\begin{align*} 1-\exp(-n^{\frac{a}{2}+\frac{\delta}{2}})\leq
x_{j^+-1}<x_{j^+}<x_{j^++1}<1.\end{align*}The proof of Lemma
\ref{lbspaceA} is completed.

From the proof of Lemma \ref{lbspaceA}, it is easy to see that each
of $[0, \exp(-n^{\frac{a}{2}+\frac{\delta}{2}})]$ and
$[1-\exp(-n^{\frac{a}{2}+\frac{\delta}{2}}), 1]$ contains a large number of
points of the minimal configuration $(x_i)_{i\in\Z}$ for $n$ large
enough.

 \vspace{2ex}

\noindent\textbf{Acknowledgement} The author would like to thank the referees for the careful reading of the paper and invaluable
comments which are very helpful in improving this paper.
This work is supported by NSFC Grant No. 11790273, 11631006.

\end{document}